\documentclass[10pt]{amsart}
 \usepackage{graphicx,enumitem,dsfont}
 \usepackage[usenames]{color}
 \usepackage[colorlinks]{hyperref}

 \newcommand{\Z}{{\mathbb{Z}} }
 
 \newcommand{\R}{\mathbb{R}}

 \newcommand{\Comment}[1]{}

 \newcommand{\eps} {\varepsilon}
 \renewcommand{\i}{\ifmmode\mathit{\mathchar"7010 }\else\char"10 \fi}
 \renewcommand{\j}{\ifmmode\mathit{\mathchar"7011 }\else\char"11 \fi}

 \newcommand{\seq}[1]{\left\{#1\right\}}

 \newcommand{\test}{\varphi}
 \newcommand{\Dx}{\Delta x}
 
 \newcommand{\Dt}{\Delta t}

 \newcommand{\udx}{u_{\Dx}}
 
 \newcommand{\norm}[1]{\left\|#1\right\|}
 \newcommand{\abs}[1]{\left|#1\right|}
 \newcommand{\abseps}[1]{\left|#1\right|_\eps}

 \newcommand{\sgn}{\mathrm{sign}}

 \newcommand{\downto}{\downarrow}
 \newcommand{\Pt}{\Pi_T}
 \newcommand{\ueta}{u^\eta}
 \newcommand{\Aeta}{A^\eta}
 \newcommand{\Dm}{D^-}
 \newcommand{\Dp}{D^+}
 \newcommand{\uDx}{u_{\Dx}}
 \newcommand{\testDx}{\test^{\Dx}}

 \DeclareMathOperator*{\sg} {sg}

 \newtheorem{definition}{Definition}[section]
 \newtheorem{theorem}{Theorem}[section]
 \newtheorem*{maincorollary*}{Main Corollary}
 \newtheorem{lemma}{Lemma}[section]
 
 \newtheorem{remark}{Remark}[section]

 \newtheorem*{maintheorem*}{Main Theorem}
 \allowdisplaybreaks
 \numberwithin{equation}{section}
 \numberwithin{figure}{section}
 \numberwithin{table}{section}
 \newcounter{asnr}
 {\ifnum\value{asnr}=0 \stepcounter{asnr} 
   \begin{enumerate}[label=\textbf{A}.\arabic{enumi}]
     \else
     \begin{enumerate}[label=\textbf{A}.\arabic{enumi},resume] \fi}
 {\end{enumerate}}
 \newcounter{defnr}
 \newenvironment{Definitions} %
 {\ifnum\value{defnr}=0 \stepcounter{defnr} 
   \begin{enumerate}[label=\textbf{D}.\arabic{enumi}]
     \else
     \begin{enumerate}[label=\textbf{D}.\arabic{enumi},resume] \fi}
 {\end{enumerate}}
 \title[Rate of convergence]{An error estimate for the finite
   difference approximation to degenerate convection - diffusion
   equations} 
 \author[K.H. Karlsen]{K. H. Karlsen} \address[Kenneth Hvistendahl
 Karlsen]{\newline 
   Centre of Mathematics for Applications (CMA) \newline University of
   Oslo\newline P.O. Box 1053, Blindern\newline N--0316 Oslo, Norway}
 \email[]{kennethk@math.uio.no}
 \author[U. Koley]{U. Koley} \address[Ujjwal Koley]{\newline Centre of
   Mathematics for Applications (CMA) 
   \newline University of Oslo\newline P.O. Box 1053, Blindern\newline
   N--0316 Oslo, Norway} \email[]{ujjwalk@cma.uio.no}
 \author[N. H. Risebro]{N. H. Risebro} \address[Nils Henrik
 Risebro]{\newline Centre of Mathematics for Applications (CMA)
   \newline University of Oslo\newline P.O. Box 1053, Blindern\newline
   N--0316 Oslo, Norway} \email[]{nilshr@math.uio.no}

\date{\today}

\begin{document}

\begin{abstract}
  We consider semi-discrete first-order finite difference schemes for
  a nonlinear degenerate convection-diffusion equations in one space
  dimension, and prove an $L^1$ error estimate. Precisely, we show
  that the $L^1_{\mathrm{loc}}$ difference between the approximate
  solution and the unique entropy solution converges at a rate
  $\mathcal{O}(\Dx^{1/11})$, where $\Dx$ is the spatial mesh size. If
  the diffusion is linear, we get the convergence rate
  $\mathcal{O}(\Dx^{1/2})$, the point being that the $\mathcal{O}$ is
  independent of the size of the diffusion.
\end{abstract}

\maketitle

\section{Introduction}
\label{sec:intro}
In this paper, we consider semi-discrete finite difference schemes for
the following Cauchy problem
\begin{equation}
  \label{eq:main}
  \begin{cases}
    u_t + f(u)_x = A(u)_{xx}, &(x,t)\in \Pt,\\
    u(0,x)=u_0(x), & x\in \R,
  \end{cases}
\end{equation}
where $\Pt=\R\times (0,T)$ with $T>0$ fixed, $u:\Pt\to \R$ is the
unknown function, $f$ the flux function, and $A$ the nonlinear
diffusion. Regarding this, the basic assumption is that $A'\ge 0$, and
thus \eqref{eq:main} is a strongly degenerate parabolic problem. The
scalar conservation law $u_t+f(u)_x=0$ is a special example of this
type of problems. Other examples occur in several applications, for
instance in porous media flow \cite{EspedalKarlsen} and in
sedimentation processes \cite{Bustosetalbook}. 

Since $A'$ may be zero, solutions are not necessary smooth and one
must consider weak solutions. These are not necessarily uniquely
determined by the initial data, and in order to get uniqueness, one
considers so-called \emph{entropy solutions}. The framework of entropy
solutions makes the initial value problem \eqref{eq:main} well posed,
for a precise statement, see Section~\ref{sec:prelim}. 
For scalar conservation laws, the entropy framework (usually called
entropy conditions) was introduced by Kru\v{z}kov \cite{Kruzkov} and
Vol'pert \cite{Volpert}, while for degenerate parabolic equations
entropy solution were first considered by Vol'pert and Hudajev
\cite{VolpertHudajev}. 
Uniqueness of entropy solutions to \eqref{eq:main} was first proved by
Carrillo \cite{Carrillo}, see also Karlsen and Risebro
\cite{KarlRise1}. 

For hyperbolic equations, the convergence analysis of difference
schemes has a long tradition, we mention only a \emph{few}
references. Finite difference schemes have been studied by Ole\u{i}nik
\cite{Oleinik}, Harten \textit{et al}.~\cite{Hartenetal}, Kuznetsov
\cite{Kuznetsov}, Crandall and Majda \cite{CrandallMajda}, Osher and
Tadmor \cite{OsherTadmor}, Cockburn and
Gripenberg.~\cite{CockburnGripenberg}, Kr\"oner and Rokyta
\cite{KronerRokyta}, Eymard \textit{et al}.~\cite{Eymardetal}, Noelle
\cite{Noelle} as well as many others.

In the last decade, there has been a growing interest in numerical
approximation of entropy solutions to degenerate parabolic
equations. Finite difference and finite volume schemes for degenerate
equations were analysed by Evje and Karlsen
\cite{EvjeKarlsen1,EvjeKarlsen2, EvjeKarlsen3, EvjeKarlsen4} (using
upwind difference schemes), Holden \textit{et
  al}.~\cite{Holdenetal1, Holdenetal2} (using operator splitting
methods),  Kurganov and Tadmor \cite{KurganovTadmor} (central
difference schemes), Bouchut \textit{et al}.~\cite{Boucutetal} (kinetic
BGK schemes), Afif and Amaziane \cite{AfifAmaziane} and Ohlberger, Gallou\"{e}t
\textit{et
  al}.~\cite{Ohlberger, Gall-1, Gall-2} (finite volume methods), Cockburn and Shu
\cite{CockburnShu} (discontinuous Galerkin methods) and Karlsen and
Risebro \cite{KRT1,KRT2} (monotone difference schemes). Many of the
above papers show that the approximate solutions converge to the
unique entropy solution as the discretization parameter vanishes. 

Despite this relatively large body of research, to the best of our
knowledge, there does not exist a result giving the \emph{convergence
  rate} of the approximate solutions to degenerate problems. For
conservation laws (very degenerate problems), the convergence rate for
monotone methods has long been known to be $\Dx^{1/2}$
\cite{Kuznetsov}, and this is also optimal for discontinuous
solutions. For non-degenerate problems, the solution operator (taking
initial data to the corresponding solution) has a strong smoothing
effect, and truncation analysis applies. Hence difference methods
produces approximations converging at the formal order of the
scheme. However, all estimates depend on  $A'$, and are not
available if $A'$ is not bounded below by some positive number
$\eta$. 

Often, the viscous regualarization 
\begin{equation}
u^{\eta}_t + f(u^\eta)_x = A(u^\eta)_{xx} + \eta u^\eta_{xx},\label{eq:viscreg}
\end{equation}
is used to model the behaviour of first order difference schemes. The
rationale behind this is that first order schemes for \eqref{eq:main}
are formally second order accurate for an equation resembling
\eqref{eq:viscreg}. If one can prove convergence, or find a
convergence rate, such that $u^\eta\to u$, then often analogous
arguments will work for appropriate difference schemes. If this
convergence has a rate, it is expected that the difference scheme will
have the same rate. In \cite{EvjeKarlsen4}, Evje and Karlsen showed that  
\begin{equation*}
  \norm{u^\eta(\cdot,t)-u(\cdot,t)}_{L^1(\R)} =
  \mathcal{O}\left(\eta^{1/2}\right).
\end{equation*}
Also, in \cite{Gall-2},  Gallou\"{e}t \textit{et al}. showed for the \emph{boundary value} problem corresponding to 
\eqref{eq:main} that
\begin{equation*}
  \norm{u^\eta(\cdot,t)-u(\cdot,t)}_{L^1(\R)} =
  \mathcal{O}\left(\eta^{1/5}\right).
\end{equation*}
To show the same rate for a difference scheme seems remarkably
difficult. The main result of this paper is that we prove a
significantly lower convergence rate for a semi-discrete difference
approximation $\uDx$,
\begin{equation*}
  \int_{-L+Mt}^{L-Mt} \abs{u(x,t)-\uDx(x,t)} \,dx \le
  \mathcal{O}\left(\Dx^{1/11}\right), 
\end{equation*}
where $M$ is a constant larger than $\abs{f'}$. 
The $\mathcal{O}$ symbol depends on $t$, $L$ and the initial data $u_0$,
but not on the discretization parameter $\Dx$. 

The rest of this paper is organized as follows. In
Section~\ref{sec:prelim} we make precise the definition of a solution
of \eqref{eq:main}, and of $\uDx$. Then we list a number of useful
properties of the (unique) weak solution and of the approximation
$\uDx$. Finally we state our main theorem. Section~\ref{sec:proof} is
devoted to its proof, while we test the practical convergence
properties of the scheme on a numerical example in 
Section~\ref{sec:numtest}.

\section{Preliminaries}\label{sec:prelim}
Independently of the smoothness of the initial data, due to the
degeneracy of the diffusion, jumps may form in the solution
$u$. Therefore we consider solutions in the weak sense, i.e.,
\begin{definition}
  Set $\Pt=(0,T)\times\R$, a function $u(t,x) \in L^{\infty}\left((0,T);L^1(\R)\right)\cap
  L^\infty(\Pt)$ is a weak solution of the initial value problem
  \eqref{eq:main} if it satisfies
  \begin{Definitions}
  \item $A(u)$ is continuous and $A(u)_x \in
    L^\infty(\Pt)$.\label{def:w1}
  \item For all test functions $\test\in
    \mathcal{D}(\Pt)$ \label{def:w2}
    \begin{equation}
      \label{eq:weaksol}
      \iint_{\Pt} u\test_t + f(u)\test_x + A(u)\test_{xx}\,dxdt = 0.
    \end{equation}
  \item The initial condition is satisfied in the $L^1$-sense
    \begin{equation*}
      \lim_{t\downto 0} \int_{\R} \abs{u(t,x)-u_0(x)}\,dx = 0. 
    \end{equation*}
  \end{Definitions}
\end{definition}
In view of the existence theory, the condition~\ref{def:w1} is
natural, and thanks to this we can replace \eqref{eq:weaksol} by
\begin{equation}
  \label{eq:weaksol2}
  \iint_{\Pt} u\test_t +\left( f(u) - A(u)_x\right) \test_x \,dxdt = 0.
\end{equation}
If $A$ is constant on a whole interval, then weak solutions are not
uniquely determined by their initial data, and one must impose an
additional entropy condition to single out the physically relevant
solution. A weak solution satisfies the entropy condition if
\begin{equation}
  \label{eq:entcond1}
  \varrho(u)_t + q(u)_x + r(u)_{xx}\le 0\ \text{in $\mathcal{D}'(\Pt)$,}
\end{equation}
for all convex, twice differentiable functions $\varrho:\R\to \R$,
where $q$ and $r$ are defined by
$$
q'(u)=\varrho'(u)f'(u), \ \text{and}\ r'(u)=\varrho'(u) A'(u).
$$
Via a standard limiting argument this implies that \eqref{eq:entcond1}
holds for the Kru\v{z}kov entropies $\varrho(u)=\abs{u-c}$ for all
constants $c$. We say that a weak solution satisfying the entropy
condition is an entropy solution.

Let the signum function be defined as
$$
\sgn(\sigma) =
\begin{cases}
  -1 & \sigma<0,\\
  0 & \sigma = 0,\\
  1 & \sigma>0,
\end{cases}
$$
and its regularized counterpart, $\sgn_\eps$, defined as
$$
\sgn_\eps(\sigma) =
\begin{cases}
  \sgn(\sigma) & \abs{\sigma}>\eps,\\
  \sin\left(\frac{\pi\sigma}{2\eps}\right) &\text{otherwise,}
\end{cases}
$$
where $\eps>0$.

We collect some useful information about entropy solutions in the
following, for a proof see \cite{EvjeKarlsen4}.
\begin{theorem}
  The unique entropy solution $u$ of \eqref{eq:main} satisfies
  \begin{equation}
    \label{eq:ws_etropy1}
    \iint_{\Pt} \abs{u-c}\test_t + \sgn(u-c)(f(u)-f(c))\test_x +
    \abs{A(u)-A(c)} \test_{xx} \,dxdt \ge 0,
  \end{equation}
  for all constants $c$ and all non-negative test functions in
  $\mathcal{D}'(\Pt)$.  Furthermore, the following limits hold, provided $A$ is strictly increasing, 
  \begin{equation}
    \label{eq:ws_limit1}
    \begin{aligned}
      \iint_{\Pt} \abs{u-c}\test_t +& \sgn(u-c)\left(f(u)-f(c) -
        A(u)_x\right) \test_x\\ & = \lim_{\eps\downto 0} \iint_{\Pt}
      \abs{A(u)_x}^2 \sgn'_{\eps}\left(A(u)-A(c)\right)\,\test \,dtdx,
    \end{aligned}
  \end{equation}
  \begin{equation}
    \label{eq:we_limit2}
    \lim_{\eps\downto 0} \iint_{\Pt} \left(f(u)-f(c)\right)A(u)_x\,
    \sgn'_{\eps} \left(A(u)-A(c)\right)\test \,dtdx = 0,
  \end{equation}
  for all non-negative test functions $\test$.
\end{theorem}
A common method to show existence of an entropy solution is to
consider the regularized problem
\begin{equation}
  \label{eq:regularized}
  \ueta_t + f\left(\ueta\right)_x = \left(A\left(\ueta\right)+\eta
    u\right)_{xx}, \ t>0, \quad \ueta(0,x)=u_0(x),
\end{equation}
where $\eta$ is some (small) positive number. This equation is not
degenerate, and has a unique smooth solution for $t>0$. The sequence
$\seq{\ueta}_{\eta>0}$ is compact in $L^1(\Pt)$, and converges to the
entropy solution. In \cite{EvjeKarlsen4}, it was established that for
$t<T$
\begin{equation}
  \label{eq:diff_u_ueta}
  \norm{u(t,\cdot)-\ueta(t,\cdot)}_{L^1(\R)} \le C \sqrt{\eta},
\end{equation}
where the constant $C$ only depends on $f$, $A$ and the initial data
$u_0$. Of course, $\ueta$ is also an entropy solution to
\eqref{eq:regularized}, where the diffusion function is given by
$$
\Aeta(u)=A(u)+\eta u.
$$

Rather than discretizing \eqref{eq:main}, we shall discretize the
regularized equation \eqref{eq:regularized}, and let $\eta$ tend to
zero in a suitable manner. Due to \eqref{eq:diff_u_ueta}, it suffices
to compare $\ueta$ and our approximate solution. To simplify our
notation, we therefore, \emph{for the moment}, assume that $A$ is
strictly increasing, with $A'(u)\ge \eta>0$.

We consider a semi-discrete approximation, where space is
discrete, but time continuous. Let $\Dx$ be some small parameter, and
set $x_j=j\Dx$ and $x_{j+1/2}=(j+1/2)\Dx$ for $j\in\Z$. Set $I_j =
(x_{j-1/2},x_{j+1/2}]$. The discrete derivatives $D^{\pm}$ are
defined by
$$
D^{\pm}\sigma_j =\pm \frac{\sigma_{j\pm 1} - \sigma_j}{\Dx}.
$$
Our scheme is defined by
\begin{equation}
  \label{eq:scheme}
  \begin{cases}
    \frac{d}{dt} u_j(t) + \Dm F_{j+1/2} = \Dm\Dp A_j, & t>0,\\
    u_j(0)=\frac{1}{\Dx}\int_{I_j} u_0(x)\,dx,
  \end{cases}
\end{equation}
for $j\in\Z$. Here $F_{i+1/2}$ is the Engquist-Osher flux and $A_j = A\left(u_j\right)$.
More precisely, for a monotone flux $f$, the generalized upwind scheme of Engquist and
Osher is defined by
\begin{align*}
F_{j+\frac{1}{2}} = F(u_j, u_{j+1}) = f^{+}(u_j) + f^{-}(u_{j+1}),
\end{align*}
where
\begin{align*}
f^{+}(u) = f(0) + \int_0^u \max \left( f'(s),0\right)\,ds, \qquad f^{-}(u) = \int_0^u \min \left( f'(s),0\right)\,ds.
\end{align*}
With this, we can rewrite our scheme \eqref{eq:scheme} as
\begin{equation}
  \label{eq:scheme_linear}
  \begin{cases}
    \frac{d}{dt} u_j(t) + \left( \Dm f^{+}(u_j) + \Dp f^{-}(u_{j})\right) = \Dm\Dp A_j, & t>0,\\
    u_j(0)=\frac{1}{\Dx}\int_{I_j} u_0(x)\,dx,
  \end{cases}
\end{equation}
for $j\in\Z$ and $A_j = A\left(u_j\right)$.

In order to define an approximation on the whole of $\Pt$, we let
$\uDx$ be the piecewise linear interpolant given by
$$
\uDx(x,t)=u_j(t) + \Dp u_j(t) \left(x-x_j\right), \quad \text{for
  $x\in [x_j,x_{j+1}]$,}
$$
and with a slight abuse of notation we define $u_j$ to be the
piecewise constant (in $x$) function
$$
u_j(x,t)=u_j(t) \quad \text{for $x\in (x_{j-1/2},x_{j+1/2}]$.}
$$
We collect some useful results regarding $\uDx$ and the entropy
solution $u$ in the next lemma.
\begin{lemma}
  \label{lem:facts} If $u$ is the unique entropy solution of
  \eqref{eq:main} and $u_j$ the function defined by the scheme
  \eqref{eq:scheme}. Then the following estimates hold:
  \begin{align}
    \norm{u(\cdot,t)}_{L^\infty(\R)} &\le
    \norm{u_0}_{L^\infty(\R)} \label{eq:linftyu}
    \\
    \abs{u(\cdot,t)}_{B.V.(\R)} &\le
    \abs{u_0}_{B.V.(\R)}\label{eq:bvu}
    \\
    \norm{f(u(\cdot,t))-A(u(\cdot,t))_x}_{L^\infty(\R)} &\le
    \norm{f(u_0)-A(u_0)_x}_{L^\infty(\R)} \label{eq:fluxlinftyu}
    \\
    \abs{f(u(\cdot,t))-A(u(\cdot,t))_x}_{B.V.(\R)} &\le
    \abs{f(u_0)-A(u_0)_x}_{B.V.(\R)} \label{eq:fluxbvu}
    \\
    \norm{u_j(t)}_{L^\infty(\R)} &\le
    \norm{u_j(0)}_{L^\infty(\R)}\label{eq:linftyudx}
    \\
    \abs{u_j(t)}_{B.V.(\R)} &\le
    \abs{u_j(0)}_{B.V.(\R)}\label{eq:bvudx}
    \\
    \norm{F_{j+1/2}(t) - \Dp A\left(u_j(t)\right)}_{L^\infty(\R)} &\le
    \norm{F_{j+1/2}(0) - \Dp
      A\left(u_j(0)\right)}_{L^\infty(\R)} \label{eq:fluxlinftyudx}
    \\
    \abs{F_{j+1/2}(t) - \Dp A\left(u_j(t)\right)}_{B.V.(\R)} &\le
    \abs{F_{j+1/2}(0) - \Dp
      A\left(u_j(0)\right)}_{B.V.(\R)}. \label{eq:fluxbvudx}
  \end{align}
  In addition the initial error is bounded by
  \begin{equation}
    \label{eq:initial_error}
    \norm{\uDx(\cdot,0)-u_0}_{L^1(\R)} \le C\Dx,
  \end{equation}
  for some constant $C$ which depends only on $\abs{u_0}_{B.V.(\R)}$.
\end{lemma}
For a proof of this, see \cite{EvjeKarlsen4}. Note that if $u_0$ and
$A(u_0)_x$ (and $u_j(0)$ and $\Dp A(u_j(0))$) are bounded
independently of $\eta$ (and $\Dx$), then \eqref{eq:fluxlinftyu} and
\eqref{eq:fluxlinftyudx} imply that
\begin{equation}
  \norm{\Dp A(u_j)}_{L ^\infty(\R)} \le C \ \text{and}\ 
  \norm{ A(u)_x}_{L^\infty(\R)} \le C,
\end{equation}
for some constant $C$ which is independent of $\Dx$ and $\eta$.

Our main result is the following
\begin{maintheorem*}
  Let $u$ be the unique entropy solution to \eqref{eq:main} and $\uDx$
  be as defined by \eqref{eq:scheme_linear}. Choose a constant 
  $$
  M>\max_{\abs{u}<\norm{u_0}_{L^\infty(\R)}} \abs{f'(u)},
  $$
  and another constant $L>MT$, where $T>0$.
  Then there exists a constant $C$, independent
  of $\Dx$, but depending on $f$, $L$, $T$ and $u_0$, such that
  \begin{equation*}
    \int_{L+Mt}^{L-Mt}\abs{u(t,x) - \uDx(t,x)}\,dx  \le C \Dx^{1/11}
    \quad\text{for $t\le T$.}
  \end{equation*}
\end{maintheorem*}
As a by-product of our method of proof we get an improved rate if the
diffusion is linear. The significance of this rate is that is
independent of the size of the diffusion, which in this case is
$\eta$. 
\begin{maincorollary*}
  Let $u$ be the unique solution to the viscous regularization
  \begin{equation*}
    u_t + f(u)_x = \eta u_{xx},\ t>0,\quad u(x,0)=u_0(x),
  \end{equation*}
  and let $\uDx$ be defined by \eqref{eq:scheme_linear} with $A(u)=\eta
  u$. Then there exists a constant $C$, independent of $\Dx$ and
  $\eta$, but depending on $f$, $L$, $T$ and $u_0$, such that
   \begin{equation*}
    \int_{L+Mt}^{L-Mt}\abs{u(t,x) - \uDx(t,x)}\,dx  \le C \Dx^{1/2}
    \quad\text{for $t\le T$.}
  \end{equation*}
\end{maincorollary*}
In the linear case this is what we expect. In fact, in
\cite{Chen:2006oy} Chen and Karlsen showed that for a \emph{linear}
flux function $f$ the expected $\eta$ independent rate of $1/2$
holds. However, the methods used in \cite{Chen:2006oy} are not easily
modified to nonlinear flux functions.

\section{Proof of the main theorem}\label{sec:proof} 
First of all, for simplicity we 
will prove the main theorem for $ f' <0 $. Note that, in that case $F_{j+1/2} = f(u_{j+1})$.
The general case treatment will be similar (see Remark ~\ref{remark:1}).
The theorem will be proved by a ``doubling of the variables''
argument, but we start not with the entropy condition
\eqref{eq:ws_limit1}, but in the argument leading up to this
condition. Set
\begin{equation*}
  \psi_\eps(u,c)=\int_c^u \sgn_\eps\left(A(z)-A(c)\right)\,dz.
\end{equation*}
This is a convex entropy for all constants $c$. Set $u=u(y,s)$ and
rewrite \eqref{eq:main} as
\begin{equation*}
  u_s + \left(f(u)-f(c)\right)_y = \left(A(u)-A(c)\right)_{yy},
\end{equation*}
and multiply this with $\psi'_\eps(u,c)\test$ where $\test$ is a test
function with compact support in $\R\times (0,T)$.  Remember that
$A'\ge \eta>0$, so $u$ is smooth, after a partial
integration, we arrive at
\begin{align*}
  &\iint_{\Pt} \psi_\eps(u,c) \test_s +
   Q_{\epsilon} (u,c) \test_y \,dyds \\
  &\qquad = \iint_{\Pt} \sgn_\eps\left(A(u)-A(c)\right)
  A(u)_y\test_y + \sgn_\eps'\left(A(u)-A(c)\right)
  \left(A(u)_y\right)^2\test \,dyds.
\end{align*}
Where we have used $Q'_{\eps}(u,c) = \psi'_\eps(u,c) f'(u)$. Although $\psi_\eps(u,c)\approx \abs{u-c}$, $\psi_\eps$ is not
symmetric in $u$ and $c$. This makes it cumbersome to work with when
doubling the variables, so we rewrite the above as
\begin{equation}
  \begin{aligned}
    \iint_{\Pt} &\abs{u-c} \test_s +
          Q_{\epsilon} (u,c) \test_y \,dyds \\
    &= \iint_{\Pt} \sgn_\eps'\left(A(u)-A(c)\right)
    \left(A(u)_y\right)^2\test + \sgn_\eps\left(A(u)-A(c)\right)
    A(u)_y\test_y
    \\
    &\qquad \hphantom{\iint_{\Pt} -} \qquad +
    \left(\abs{u-c}-\psi_\eps(u,c)\right)\test_s \,dyds.
  \end{aligned}\label{eq:upre_entropy}
\end{equation}
In the doubling of variables argument we choose $c=\uDx(x,t)$, and a
test function $\test(x,y,t,s)$. Integrating the above over $(x,t)\in
\Pt$ after an integration by parts, we end up with
\begin{multline}
  \label{eq:udoubled}
  \int_{\Pt^2} \abs{u-\udx}\test_s +
  \sgn_{\eps} (A(u) -A(\uDx)) (f(u)-f(\uDx)) \test_y \,dX
  \\
  =\int_{\Pt^2}\Bigl[\sgn'_\eps\left(A(u)-A\left(\uDx\right)\right)
  \left(\left(A(u)_y\right)^2-A(u)_yA\left(\uDx\right)_x\right) \test
  \\
  -
  \abseps{A(u)-A\left(\uDx\right)}\left(\test_{yy}+\test_{xy}\right)
  + \left(\psi_\eps\left(u,\uDx\right)-\abs{u-\uDx}\right)\test_s
\\
+
  \left(\int_{\uDx}^u \frac{d}{dz}(\sgn_{\eps} (A(z) -A(\uDx))) (f(z)-f(\uDx)) \,dz \right) \test_y 
  \Bigr] \,dX,
\end{multline}
where $dX=dydsdxdt$ and $\abseps{a}=\int^a_0 \sgn_\eps(z)dz$. Here we
have used that 
\begin{align*}
  0&=\int_{\Pt^2}  \left(\sgn_\eps\left(A(u)-A(\uDx)\right) A(u)_y
    \test \right)_x \, dX \\
   &= \int_{\Pt^2} \sgn_\eps'\left(A(u)-A(\uDx)\right) A(u)_y
   A(\uDx)_x \test
   -  \abseps{A(u)-A\left(\uDx\right)} \test_{xy} \, dX,
\end{align*}
and that 
\begin{align*}
  \int_{\Pt^2} \sgn_\eps\left(A(u)-A(\uDx)\right) A(u)_y \test_y \,dX
  &= \int_{\Pt^2} \left( \abseps{A(u)-A\left(\uDx\right)}\right)_y
  \test_y \,dX\\
  &= - \int_{\Pt^2}  \abseps{A(u)-A\left(\uDx\right)}
  \test_{yy} \,dX.
\end{align*}
Also 
 \begin{align*}
Q_{\eps}(u,\uDx) & = \int_{\uDx}^u \sgn_{\eps} (A(z) -A(\uDx)) \frac{d}{dz}(f(z)-f(\uDx)) \,dz \\
& = - \int_{\uDx}^u \frac{d}{dz}(\sgn_{\eps} (A(z) -A(\uDx))) (f(z)-f(\uDx)) \,dz \\
&\qquad \qquad + \sgn_{\eps} (A(u) -A(\uDx)) (f(u)-f(\uDx)).
\end{align*}
The next goal is to obtain an analogous estimate for the difference
approximation $\uDx$. Set $\test_j(t)=\test(x_j,t)$ and multiply the
scheme \eqref{eq:scheme} with $\psi_\eps'(u_j,c)\test_j$ and do a
summation by parts to get
\begin{align*}
  &\sum_j \Bigl[\psi\left(u_j,c\right)_t \test_j  + \test_j \Dp Q_{\epsilon}(u_j,c) \Bigr]\\
  &= -\sum_j \psi_\eps'\left(u_j,c\right) \Dp
  A\left(u_j\right)\Dp\test_j + \Dp\psi_{\eps}'\left(u_j,c\right) \Dp
  A\left(u_j\right) \test_{j+1}\\
& \qquad \qquad + \sum_j  \frac{\test_j}{\Dx} \int_{u_j}^{u_{j+1}} \psi_\eps{''}\left(s,c\right)(f(u_{j+1}) -f(s))\,ds.
\end{align*}
Where we have used the following result:
\begin{align*}
0 \ge \int_{u_j}^{u_{j+1}} \psi_{\eps}^{''} (s,c) [f(u_{j+1}) -f(s)]\,ds = \int_{u_j}^{u_{j+1}}  \psi_{\eps}^{'} (s,c) & f'(s)\,ds \\
& - \psi_{\eps}'(u_j) [f(u_{j+1}) -f(u_j)]. 
\end{align*}
Set $\testDx=\test_j$ for $x\in (x_{j-1/2},x_{j+1/2}]$, integrating
the above for $t\in [0,T]$ and multiplying with $\Dx$, we obtain
\begin{align*}
  \iint_{\Pt} &\abs{\uDx-c} \testDx_t +
  \sgn_\eps\left(A\left(u_j\right)-A(c)\right)\left(f(u_j)-
    f(c)\right) \Dm \testDx \,dxdt
  \\
  &\qquad +\iint_{\Pt}
  \left(\psi_\eps\left(u_j,c\right)-\abs{\uDx-c}\right)\testDx_t
  \,dxdt \\
  &\qquad -\iint_{\Pt}
  \sgn_{\eps}\left(A\left(u_j\right)-A(c)\right)\Dp
  A\left(u_j\right)\Dp \testDx \,dxdt
  \\
  & \ge \iint_{\Pt}
  \sgn_\eps'\left(A\left(\theta_{j+1/2}\right)-A(c)\right) \left[\Dp
    A\left(u_j\right)\right]^2 \test_{j+1} \,dxdt
  \\
  &\hphantom{\iint_{\Pt}}+ \iint_{\Pt} 
   \left(\int_{c}^{u_j} \frac{d}{dz}(\sgn_{\eps} (A(z) -A(c))) (f(z)-f(c)) \,dz \right) \Dm \testDx  \,dxdt.
\end{align*}
This can be rewritten
\begin{align*}
  \iint_{\Pt} & \abs{\uDx-c}\test_t +
  \sgn_\eps\left(A(u_j)-A(c)\right)
  \left(f\left(u_j\right)-f(c)\right)\test_x \,dxdt \\
  & \ge  \iint_{\Pt} \sgn_{\eps}\left(A\left(u_j\right)-A(c)\right)\Dp
  A\left(u_j\right)\test_x \,dxdt\\
  &\quad + \iint_{\Pt}
  \sgn_\eps'\left(A\left(\theta_{j+1/2}\right)-A(c)\right)
  \left[\Dp A\left(u_j\right)\right]^2 \test \,dxdt \\
  &\quad + \iint_{\Pt}
  \sgn_\eps'\left(A\left(\theta_{j+1/2}\right)-A(c)\right)
  \left[\Dp A\left(u_j\right)\right]^2 \left(\test_{j+1}-\test\right) \,dxdt \\
  &\quad + \iint_{\Pt}
  \sgn_{\eps}\left(A\left(u_j\right)-A(c)\right)\Dp
  A\left(u_j\right)\left(\Dp \testDx-\test_x\right) \,dxdt\\
  &\quad + \iint_{\Pt}
  \left(\abs{\uDx-c}-\psi_\eps\left(u_j,c\right)\right) \testDx_t +
  \abs{\uDx-c} \left(\test_t-\testDx_t\right) \, dxdt\\
  &\quad + \iint_{\Pt}
  \sgn_\eps\left(A(u_j)-A(c)\right)\left(f(u_j)-f(c)\right)
  \left(\test_x -\Dm\testDx_j\right)\,dxdt\\
  &\quad + \iint_{\Pt} 
   \left(\int_{c}^{u_j} \frac{d}{dz}(\sgn_{\eps} (A(z) -A(c))) (f(z)-f(c)) \,dz \right) \Dm \testDx  \,dxdt.
\end{align*}
In order to make this more compatible with the corresponding equality
for the exact solution \eqref{eq:upre_entropy}, we rewrite again to
get (recall that $u_j$ denotes the piecewise constant function taking
the value $u_j(t)$ in the cell $(x_{j-1/2},x_{j+1/2}]$),
\begin{align}
  \iint_{\Pt} & \abs{\uDx-c}\test_t +
  \sgn_\eps\left(A(\uDx)-A(c)\right)
  \left(f\left(\uDx\right)-f(c)\right)\test_x \,dxdt\label{eq:udx-ent-left} \\
  & \ge  \iint_{\Pt} \sgn_{\eps}\left(A\left(\uDx\right)-A(c)\right)
  A\left(\uDx\right)_x\test_x \,dxdt\label{eq:udx-ent-r1}\\
  &\quad + \iint_{\Pt} \sgn_\eps'\left(A\left(\uDx\right)-A(c)\right)
  \left[A\left(\uDx\right)_x\right]^2 \test
  \,dxdt \label{eq:ude-ent-r2}\\
  &\quad + \iint_{\Pt} \left[
    \sgn_\eps\left(A\left(\uDx\right)-A(c)\right)-
    \sgn_\eps\left(A(u_j)-A(c)\right) \right]\notag\\
  &\hphantom{\quad + \iint_{\Pt} } \qquad\times
  \left(f\left(\uDx\right)-f(c)\right)
  \test_x\,dxdt \label{eq:udx-ent-r3} \\
  &\quad + \iint_{\Pt} \sgn_\eps\left(A\left(u_j)-A(c)\right)\right)
  \left(f\left(\uDx\right)-f\left(u_j\right)\right) \test_x\,dxdt
  \label{eq:udx-ent-r4}\\
  &\quad + \iint_{\Pt}
  \left[\sgn_\eps'\left(A\left(\theta_{j+1/2}\right)-A(c)\right)
    -\sgn_\eps'\left(A\left(\uDx\right)-A(c)\right)\right]
  \notag\\
  &\hphantom{\quad + \iint_{\Pt} } \qquad\times\left[\Dp
    A\left(u_j\right)\right]^2 \test
  \,dxdt \label{eq:udx-ent-r5}\\
  &\quad + \iint_{\Pt} \sgn'_\eps\left(A\left(\uDx\right)-A(c)\right)
  \left[ \left(\Dp A\left(u_j\right)\right)^2 -
    \left(A\left(\uDx\right)_x\right)^2
  \right] \test\, dxdt \label{eq:udx-ent-r6}\\
  &\quad + \iint_{\Pt}
  \sgn_\eps'\left(A\left(\theta_{j+1/2}\right)-A(c)\right) \left[\Dp
    A\left(u_j\right)\right]^2 \left(\test_{j+1}-\test\right) \,dxdt
  \label{eq:udx-ent-r7}\\
  &\quad + \iint_{\Pt} \left[
    \sgn_\eps\left(A\left(\uDx\right)-A(c)\right)-
    \sgn_\eps\left(A(u_j)-A(c)\right) \right]\notag\\
  &\hphantom{\quad + \iint_{\Pt} } \qquad\times
   \Dp A\left(u_j\right)\test_x\,dxdt \label{eq:udx-ent-r71}\\
  &\quad + \iint_{\Pt} \sgn_\eps\left(A\left(\uDx\right)-A(c)\right)
  \left[ \Dp A\left(u_j\right) - A\left(\uDx\right)_x \right] \test_x\, dxdt \label{eq:udx-ent-r72}\\
  &\quad + \iint_{\Pt}
  \sgn_{\eps}\left(A\left(u_j\right)-A(c)\right)\Dp
  A\left(u_j\right)\left(\Dp \testDx-\test_x\right) \,dxdt
  \label{eq:udx-ent-r8}\\
  &\quad + \iint_{\Pt}
  \left(\abs{\uDx-c}-\psi_\eps\left(u_j,c\right)\right) \testDx_t +
  \abs{\uDx-c} \left(\test_t-\testDx_t\right) \, dxdt
  \label{eq:udx-ent-r9}\\
  &\quad + \iint_{\Pt}
  \sgn_\eps\left(A(u_j)-A(c)\right)\left(f(u_j)-f(c)\right)
  \left(\test_x -\Dm\testDx_j\right)\,dxdt
  \label{eq:udx-ent-r10}\\
  &\quad + \iint_{\Pt} 
   \left(\int_{c}^{u_j} \frac{d}{dz}(\sgn_{\eps} (A(z) -A(c))) (f(z)-f(c)) \,dz \right) \Dm \testDx  \,dxdt.  \label{eq:udx-ent-r11}
\end{align}
\begin{remark}
\label{remark:1}
In the general case, we can write $F_{j+\frac{1}{2}} = f^{+}(u_j) + f^{-}(u_{j+1})$. Note that $(f^{+})' \ge 0$,
$(f^{-})' \le 0$, and both are Lipschitz continuous. Then if we multiply \eqref{eq:scheme_linear} by $\psi_\eps'(u_j,c)$,
numerical flux part can be written as $ \Dm Q_{\epsilon}^{+} + \Dp Q_{\epsilon}^{-} + \text{``terms having the right sign''}$. Here $(Q^{\pm}_{\eps})'(u,c) = \psi'_\eps(u,c) (f^{\pm})'(u)$, and $Q_{\eps}^{+} + Q_{\eps}^{-}= Q_{\eps}$.
\end{remark}
Let now
\begin{align*}
  \mathcal{R}_1(c)&=\text{\eqref{eq:udx-ent-r3}}+
  \text{\eqref{eq:udx-ent-r4}}+\text{\eqref{eq:udx-ent-r5}}+
  \text{\eqref{eq:udx-ent-r6}},\\
  \mathcal{R}_2(c)&=\text{\eqref{eq:udx-ent-r7}}+
  \eqref{eq:udx-ent-r71}+ \eqref{eq:udx-ent-r72}+\cdots+
  \eqref{eq:udx-ent-r11}.
\end{align*}
Now we choose $c=u(y,s)$ and the same test function as before, and
integrate the result over $(y,s)\in \Pt$. The result reads
\begin{multline}
  \label{eq:udxdoubled} \int_{\Pt^2} \abs{\uDx-u} \test_t +
  \sgn_\eps\left(A\left(\uDx\right)-A(u)\right)
  \left(f\left(\uDx\right)-f\left(u\right)\right) \test_x \,dX\\
  \ge  \int_{\Pt^2} \sgn'_\eps\left(A\left(\uDx\right)-A(u)\right)
  \left(\left(A\left(\uDx\right)_x\right)^2-A\left(\uDx\right)_x
    A(u)_y\right) \test\, dX\\
  -\int_{\Pt^2} \abseps{A\left(\uDx\right)-A(u)}
  \left(\test_{xx}+\test_{xy}\right)\,dX\\
  + \iint_{\Pt} \mathcal{R}_1(u) + \mathcal{R}_2(u)\, dyds.
\end{multline}
Adding this and \eqref{eq:udoubled}, we get
\begin{align}
  \int_{\Pt^2}& \Bigl[
  \abs{\uDx-u}\left(\test_t+\test_s\right) )\label{eq:estimate-start}
  \\
  &\hphantom{ \int_{\Pt^2} \Bigl[} \quad 
  + \sgn_\eps\left(A\left(\uDx\right)-A(u)\right)
  \left(f\left(\uDx\right)-f(u)\right)
  \left(\test_x+\test_y\right)\notag\\
  & \hphantom{\int_{\Pt^2} \Bigl[ } 
  \qquad  +
  \abseps{A(\uDx)-A(u)}\left(\test_{xx}+2\test_{xy}+\test_{yy}\right)
  \Bigr]\,dX
  \notag
  \\
  & \ge  \int_{\Pt^2} \sgn'_\eps\left(A\left(\uDx\right)-A(u)\right)
  \left(A\left(\uDx\right)_x-A(u)_y\right)^2 \test \,dX \notag 
  \\
  &\quad + \int_{\Pt^2}
  \left(\psi_\eps\left(u,\uDx\right)-\abs{u-\uDx}\right)\test_s
  \,dX\notag\\
  &\quad + \int_{\Pt^2}
  \left(\int_{\uDx}^u \frac{d}{dz}(\sgn_{\eps} (A(z) -A(\uDx))) (f(z)-f(\uDx)) \,dz \right) \test_y \,dX \notag \\
  &\qquad + \iint_{\Pt} \mathcal{R}_1(u) + \mathcal{R}_2(u)\, dyds
  \notag\\
  &\ge \int_{\Pt^2}
  \left(\int_{\uDx}^u \frac{d}{dz}(\sgn_{\eps} (A(z) -A(\uDx))) (f(z)-f(\uDx)) \,dz \right) \test_y \,dX \notag \\
  &\quad + \int_{\Pt^2}
  \left(\psi_\eps\left(u,\uDx\right)-\abs{u-\uDx}\right)\test_s
  \,dX\notag\\
   \qquad \qquad &+ \iint_{\Pt} \mathcal{R}_1(u) + \mathcal{R}_2(u)\, dyds
  \notag .\\
 &=: \int_{\Pt^2} \mathcal{Q}_1 + \mathcal{Q}_2 \, dX +  \iint_{\Pt}
  \mathcal{R}_1(u) + \mathcal{R}_2(u)\, dyds \notag
\end{align}

Now we are going to specify a nonnegative test function $ \test =
\test(t,x,s,y)$ defined in $\Pt\times\Pt$. Let $\omega \in
C_{0}^{\infty} (\R)$ be a function satisfying
\begin{equation*}
  \mathrm{supp}(\omega) \subset [-1,1], \qquad
  \omega(\sigma) \ge 0 , \qquad \int_{\R} \omega(\sigma)\, d\sigma = 1,
\end{equation*}
and define $\omega_r(x)=\omega(x/r)/r$.  Furthermore, let $h(z)$ be
defined as
\begin{equation*}
  h(z)=
  \begin{cases}
    0, &z<-1,\\ z+1 & z\in [-1,0],\\ 1 &z>0.
  \end{cases}
\end{equation*}
and set $h_\alpha(z)=h(\alpha z)$. 
Let $\nu<\tau$ be two numbers
in $(0,T)$, for any $\alpha> 0$ define
\begin{equation*}
  \begin{gathered}
    H_{\alpha}(t)= \int_{-\infty}^{t}
    \omega_{\alpha}(\xi)\, d\xi, \\
    \begin{aligned}
      \Psi(x,t) &= \left(H_{\alpha_0} (t - \nu) - H_{\alpha_0} (t -
        \tau)\right) \left(h_\alpha(x-L_l(t)) -
        h_\alpha(x-L_r(t)-\frac{1}{\alpha})\right)\\
      &=:{\chi^{\alpha_0}_{(\nu,\tau)}(t)}\,{\chi^\alpha_{(L_l,L_r)}(x,t)}
    \end{aligned}
  \end{gathered}
\end{equation*}
where the lines $L_{l,r}$ are given by
$$
L_l(t)=-L+Mt,\ L_r(t)=L-Mt
$$
where $M$ and $L$ are positive numbers, $M$ will be specified below.
With $0<r< \min\seq{\nu, T - \tau}$ and
$\alpha_0\in(0,\min\seq{\nu-r,T-\tau-r})$ we set
\begin{equation}
  \label{eq:testfn}
  \test(x,t,y,s)=\Psi(x,t)\,\omega_r(x-y)\,\omega_{r_0}(t-s).
\end{equation}
We note that $\phi$ has compact support and also that we have,
\begin{align*}
  \test_t + \test_s &= \Psi_t(x,t)\,
  \omega_r(x-y)\,\omega_{r_0}(t-s),\\
  \test_x+\test_y&=\Psi_x(x,t)\, \omega_r(x-y)\,\omega_{r_0}(t-s),\\
  \test_{xx}+2\test_{xy}+\test_{yy}&
=\Psi_{xx}(x,t)\, \omega_r(x-y)\,\omega_{r_0}(t-s).
\end{align*}
For the record, we note that
\begin{equation}
  \label{eq:Psider}
  \begin{aligned}
    \Psi_t(x,t)&=-\chi^{\alpha_0}_{(\nu,\tau)}(t) M
    \left(h_\alpha'(x-L_l(t)) + h_\alpha'(x-L_r(t)-\frac{1}{\alpha})\right) \\
    & \qquad + 
    \left(\omega_{\alpha_0}(t-\nu) - \omega_{\alpha_0}(t-\tau)\right)
    \chi^\alpha_{(L_l,L_r)}(x,t),\\
    \Psi_x(x,t)  &= \chi^{\alpha_0}_{(\nu,\tau)}(t) 
    \left(h_\alpha'(x-L_l(t)) - h_\alpha'(x-L_r(t)-\frac{1}{\alpha})\right),
    \\
    \Psi_{xx}(x,t)&=
    \chi_{(\nu,\tau)}^{\alpha_0}(t)\left(h_\alpha^{''}(x-L_l(t))-h_\alpha^{''}(x-L_r(t)-\frac{1}{\alpha})\right).
  \end{aligned}
\end{equation}
We shall let all the ``small parameters'' $ \alpha$, $\alpha_0$,  $r$,
$r_0$, $\eps$ and $\Dx$ be sufficiently small, but fixed. The goal
of our manipulations is to obtain an inequality where the difference
between $\uDx$ and $u$ is bounded by some combination of all these
parameters.

We shall repeatedly use the fact that
\begin{equation}
  \label{eq:bvuseful}
  \begin{gathered}
    \iint \abs{v(x,t)-v(y,t)}\omega_r(x-y)\,dxdy \le C r \\
    \text{and}\ \iint \abs{v(x,s)-v(x,t)}\omega_{r_0}(t-s)\,dxds \le C
    r_0,
  \end{gathered}
\end{equation}
for $v=u$, $v=\uDx$, $v=f(u)$, $v=A(u)_x$ or $v=A(\uDx)_x$. These
estimates follow from the basic bounds 
in Lemma~\ref{lem:facts}.
Starting the first term on the left of \eqref{eq:estimate-start}, we
write
\begin{align*}
  \int_{\Pt^2} \abs{\uDx-u} & \left(\test_s  +\test_t\right)\,dX \le
  \underbrace{\int_{\Pt} \abs{\uDx(x,t)-u(x,t)} \Psi_t \,dxdt}_\delta \\
  &\quad + \underbrace{\int_{\Pt}\int_{\R} \abs{u(x,t)-u(x,s)}\abs{\Psi_t(x,t)}
  \omega_{r_0}(t-s) \,dsdxdt}_{\beta} \\ 
  &\quad + \underbrace{\int_{\Pt^2}  \abs{u(x,s)-u(y,s)} \abs{\Psi_t(x,t)}
  \omega_{r_0}(t-s)\,\omega_r(x-y)\, dX}_\gamma . 
\end{align*}
To estimate $\beta$ and $\gamma$ we use 
\begin{equation*}
  \abs{\Psi_t} \le
  \left(\omega_{\alpha_0}(t-\nu)+\omega_{\alpha_0}(t-\tau)\right) + 
  M\left(h_\alpha'(x-L_l(t)) + h_\alpha'(x-L_r(t)-\frac{1}{\alpha}) \right),
\end{equation*}
and that 
\begin{equation*}
  \int \omega_{\alpha_0}(t-\nu)\,dt \le 1\ \text{and}\ 
  \abs{\int M h'_{\alpha}(x-L_{l,r}(t))\, dt} \le C.
\end{equation*}
A typical term in $\beta$ reads
\begin{align*}
  \iiint &\abs{u(x,t)-u(x,s)} \omega_{\alpha_0}(t-\nu)
  \omega_{r_0}(t-s)\,dxdsdt\\ &\le C\iint \abs{t-s} \omega_{r_0}(t-s)
  \,\omega_{\alpha_0}(t-\nu)\,dsdt\\
  &\le C r_0,
\end{align*}
Hence 
\begin{equation*}
  \beta\le C r_0.
\end{equation*}
Similarly a typical term in $\gamma$ can be estimated
\begin{equation*}
  \iint \iint \abs{u(x,s)-u(y,s)} \omega_r(x-y) \omega_{r_0}(t-s)
  \omega_{\alpha_0}(t-\nu) \, dxdy dtds \le C r.
\end{equation*}
Thus we find that 
\begin{equation}
  \beta+\gamma \le C\left(r_0+r\right).\label{eq:betagamma}
\end{equation}
To continue the estimate with the first term on the left of
\eqref{eq:estimate-start}, we split $\delta$ as follows
\begin{align*}
  \delta&= \underbrace{-\iint_{\Pt} \chi_{(\nu,\tau)}^{\alpha_0}(t) M
    \left(h'_\alpha(x-L_l(t)) + h'_\alpha(x-L_r(t)-\frac{1}{\alpha})\right)
    \abs{\uDx(x,t)-u(x,t)} \, dxdt}_{\delta_1} \\
  &\quad + \underbrace{\iint_{\Pt} \chi^\alpha_{(L_l,L_r)}(x,t)
  \abs{\uDx(x,t)-u(x,t)}
  \left(\omega_{\alpha_0}(t-\nu)-\omega_{\alpha_0}(t-\tau)\right) \,dxdt}_{\delta_2}.
\end{align*}
The term $\delta_1$ will be balanced against the first order
derivative term on the
left hand side of \eqref{eq:estimate-start}.  To estimate $\delta_2$
we set $e(x,t)=\abs{\uDx(x,t)-u(x,t)}$ and proceed as follows
\begin{align*}
  \iint_{\Pt} &\chi^\alpha_{(L_l,L_r)}(x,t) e(x,t)
  \,\omega_{\alpha_0}(t-\nu)\,dxdt\\ &\le  
  \int \chi^\alpha_{(L_l,L_r)}(x,\nu) e(x,\nu)\,dx 
  \\ & \qquad \qquad + 
  \iint_{\Pt} \chi^\alpha_{(L_l,L_r)}(x,t) \abs{e(x,t)-e(x,\nu)} 
  \omega_{\alpha_0}(t-\nu) \,dxdt \\
  &\le \int \chi^\alpha_{(L_l,L_r)}(x,\nu) e(x,\nu)\,dx + C\alpha_0
  \\   \intertext{and similarly}
  \iint_{\Pt} &\chi^\alpha_{(L_l,L_r)}(x,t) e(x,t)
  \,\omega_{\alpha_0}(t-\tau)\,dxdt \ge 
   \int \chi^\alpha_{(L_l,L_r)}(x,\tau) e(x,\tau)\,dx - C\alpha_0.
\end{align*}
Using this we get the estimate
\begin{equation}
  \label{eq:delta2est}
  \begin{aligned}
    \delta_2 &\le \int \chi^\alpha_{(L_l,L_r)}(x,\nu) \abs{\uDx(x,\nu)-u(x,\nu)}
    \,dx \\ & \qquad - \int \chi^\alpha_{(L_l,L_r)}(x,\tau)
    \abs{\uDx(x,\tau)-u(x,\tau)} \,dx +
    C\alpha_0.
  \end{aligned}
\end{equation}
Now we rewrite the ``first derivative term'' on the left hand side of
\eqref{eq:estimate-start}. Doing this, we get
\begin{align*}
  \int_{\Pt^2} & \sgn_\eps\left(A\left(\uDx\right)-A(u)\right)
    \left(f\left(\uDx\right)-f(u)\right)
    \left(\test_x+\test_y\right) \,dX \\ &= 
    \int_{\Pt^2} \sg(x,y,t,s)\left(f(\uDx(x,t))-f(u(x,t))\right)
    \Psi_x(x,t)\,\omega_r(x-y)\, \omega_{r_0}(t-s)\,dX\\
    &\quad + 
    \int_{\Pt^2} \sg(x,y,t,s)\left(f(u(x,t))-f(u(y,s))\right) \Psi_x(x,t)\,
    \omega_r(x-y)\,\omega_{r_0}(t-s)\,dX \\
    &=:\delta_3 + \delta_4,
\end{align*}
where we have set $\sg(x,y,t,s)=\sgn_\eps(A(\uDx(x,t))-A(u(y,s)))$.
We  proceed as follows
\begin{align*}
  \abs{\delta_4} \le \int_{\Pt^2} &\abs{f(u(x,t))-f(u(y,s))}
  \chi^{\alpha_0}_{(\nu,\tau)}(t) \\
  &\qquad \omega_{r_0}(t-s)\,\omega_r(x-y)
  \left(
    h'_\alpha(x-L_l(t)) + h'_\alpha(x-L_r(t)-\frac{1}{\alpha})\right) \,dX.
\end{align*}
Each of these two terms are estimated using \eqref{eq:bvuseful} by
\begin{align*}
  \int_{\Pt^2} &\abs{f(u(x,t))-f(u(y,s))}
  \chi^{\alpha_0}_{(\nu,\tau)}(t) \,\omega_{r_0}(t-s)\,\omega_r(x-y)
  \,
  h'_\alpha(x-L_l(t))  \,dX\\
  &\le 
  \int_{\Pt^2} \abs{f(u(x,t))-f(u(x,s))}
  \chi^{\alpha_0}_{(\nu,\tau)}(t) \,\omega_{r_0}(t-s)\,\omega_r(x-y)
  h'_\alpha(x-L_l(t))\,dX \\
  &\quad + 
  \int_{\Pt^2} \abs{f(u(x,s))-f(u(y,s))}
  \chi^{\alpha_0}_{(\nu,\tau)}(t) \,\omega_{r_0}(t-s)\,\omega_r(x-y)
  h'_\alpha(x-L_l(t))\,dX \\
  &\le C\left(r_0+r\right),
\end{align*}
and thus 
\begin{equation}
  \label{eq:delta4bnd}
  \abs{\delta_4} \le C\left(r_0+r\right).
\end{equation}
The terms $\delta_1+\delta_3$ is bounded by choosing $M$ sufficiently
large (all functions are functions of $(x,t)$).
\begin{align*}
  \delta_1&+\delta_3 \\&= 
  \iint_{\Pt} \chi^{\alpha_0}_{(\nu,\tau)}(t) 
   h'_\alpha(x-L_l(t)) \left(-M\abs{\uDx-u} + \sg\left(f(\uDx)-f(u)\right)
    \right) \,dxdt \\
 &\quad + 
 \iint_{\Pt} \chi^{\alpha_0}_{(\nu,\tau)}(t) 
   h'_\alpha(x-L_r(t)-\frac{1}{\alpha}) \left(-M\abs{\uDx-u} -
     \sg\left(f(\uDx)-f(u)\right) \right) \,dxdt. 
\end{align*}
Choosing $M$ larger than the Lipschitz norm of $f$ implies that 
\begin{equation}
  \label{eq:dl1dl2bnd}
  \delta_1+\delta_3\le 0.
\end{equation}
Turning to the ``second derivative term'' \eqref{eq:estimate-start},
we get
\begin{align*}
  &\int_{\Pt^2}  \abseps{A(\uDx)-A(u)} \Psi_{xx}\,
  \omega_{r_0}(t-s)\,\omega_r(x-y)\,dX \\
  &\,\le \int_{\Pt^2} \abs{A(\uDx)-A(u)} 
  \left(h_\alpha''(x-L_l(t)) + h_\alpha''(x-L_r(t)-\frac{1}{\alpha})\right) \\
  & \qquad \qquad \qquad \qquad \qquad \qquad \omega_{r_0}(t-s)\,\omega_r(x-y)\,dX\\
  &\,\le C\alpha \int_{\Pt^2}  \left(\abs{\uDx(x,t)}+\abs{u(y,s)}\right) 
  \,\omega_{r_0}(t-s)\,\omega_r(x-y)\,dX\\
  &\,\le C\alpha,
\end{align*}
for some constant $C$ which is independent of the small parameters. 

Collecting \eqref{eq:betagamma}, \eqref{eq:delta2est},
\eqref{eq:delta4bnd}, \eqref{eq:dl1dl2bnd} and the above inequality we
see that
\begin{equation}
  \label{eq:estimate-middle}
  \begin{aligned}
    \int \chi^\alpha_{(L_l,L_r)}&(x,\tau) \abs{\uDx(x,\tau)-u(x,\tau)}
    \,dx \\
    &\le \int \chi^\alpha_{(L_l,L_r)}(x,\nu)
    \abs{\uDx(x,\nu)-u(x,\nu)} \,dx\\
    &\qquad + C\left(r_0+ r +
      \alpha_0 + \alpha\right)\\
    &\qquad + \biggl|{\int_{\Pt^2} \mathcal{Q}_1 + \mathcal{Q}_2 \, dX +
      \iint_{\Pt} \mathcal{R}_1(u) + \mathcal{R}_2(u)\, dyds }\biggr|.
  \end{aligned}
\end{equation}
In order to estimate the integral involving $\mathcal{Q}_2$ we first
observe that since $A'(u)\ge \eta$, we get
\begin{align*}
  \abs{\psi_\eps(a,b)-\abs{a-b}}&\le
  \int_b^a \abs{\sgn_\eps(A(z)-A(b)) - \sgn(A(z)-A(b))}\,dz\\
  &\le \frac{1}{\eta} \int_{A(b)}^{A(a)}
  \chi_{\abs{\beta}<\eps}(\beta) \abs{\sin\left(\frac{\pi\beta}{2\eps}\right)}\,d\beta\\
  &\le C\frac{\eps}{\eta}.
\end{align*}
Using this we find that 
\begin{align}
  \abs{\int_{\Pt^2}\mathcal{Q}_2} \,dX &\le C\frac{\eps}{r_0\eta} \int_0^T
  \int_{-L+(M/\alpha)t}^{L-(M/\alpha)t} dxdt\notag\\
  &\le C\frac{\eps}{r_0\eta}\left(L+\frac{1}{\alpha}\right),\label{eq:q2est}
\end{align}
where $C$ is independent of the small parameters. 
To estimate the integral of $\mathcal{Q}_1$ we use the Lipschitz 
continuity of $f\left(A^{-1} \right)$. Note that $f\left(A^{-1} \right)$
is Lipschitz continuous with Lipschitz constant $\frac{M}{\eta}$, where
$M$ is the Lipschitz constant for $f$. 

\begin{align}
&\int_{\uDx}^u  \abs{\frac{d}{dz}(\sgn_\eps\left(A(z)-A\left(\uDx\right)\right))}
  \abs{\left(f(z)-f\left(\uDx\right)\right)}  \,dz \notag\\
& = \int_{A(\uDx)}^{A(u)} \sgn_\eps' \left( r -A\left(\uDx\right)\right)
  \abs{\left(f \left( A^{-1}(r) \right)-f\left(\uDx\right)\right)}  \,dr \notag\\
& = \frac{1}{\epsilon} \int_{\min\left( A(u), A(\uDx) -\epsilon \right)}^{\min\left( A(u), A(\uDx) +\epsilon \right)}
 \abs{\left(f \left( A^{-1}(r) \right)-f\left(\uDx\right)\right)}  \,dr \notag\\
& \le \frac{\eps^2 M}{\eta \epsilon}
\end{align}
Using this we find that
\begin{align}
  \abs{\int_{\Pt^2} \mathcal{Q}_1\,dX} &\le 
  \frac{M\eps}{\eta} \int_{\Pt^2}  \test_y \,dX\notag\\
  &\le C \frac{\eps}{\eta r}\left(L+\frac{1}{\alpha}\right).\label{eq:q1est}
\end{align}
Therefore
\begin{equation}
  \label{eq:qsum}
  \abs{\int_{\Pt^2} \mathcal{Q}_1 + \mathcal{Q}_2 \,dX } \le 
  \frac{C\eps}{\eta}\left(L+\frac{L}{r_0} + \frac{1}{r \alpha}\right).
\end{equation}
Now we claim that 
\begin{equation}
  \label{eq:Rclaim}
  \begin{aligned}
    \abs{ \iint_{\Pt} \mathcal{R}_1(u) + \mathcal{R}_2(u)\, dyds } &\le
    C\Bigl[\Dx\left(\frac{1}{\eps^2} + \frac{1}{\eps\eta^2} +
        \frac{1}{r^2} + \frac{1}{r_0} + \frac{1}{\eps r} + \frac{1}{r} \right) \\
      &\hphantom{\le
    C\Bigl[}\quad +
      \eps\left(\frac{1}{\eta} + 
        \frac{1}{r_0\eta} + \frac{1}{r_0\eta\alpha} +
        \frac{1}{\eta\alpha} \right) \Bigr], 
  \end{aligned}
\end{equation}
where $C$ depends on (among other things) $L$ and $T$, but not on
the small parameters $\alpha_0$, $\alpha$, $r_0$, $r$, $\eta$ or $\eps$.
The proof of this claim is a tedious computation of all the terms. We
start by the one from \eqref{eq:udx-ent-r3}.
\begin{align*}
  \iint_{\Pt} \abs{\eqref{eq:udx-ent-r3}}\,\bigm|_{c=u}\,dyds&\le 
  \int_{\Pt^2}
  \abs{\sgn_\eps\left(A(\uDx)-A(u)\right)-\sgn_\eps\left(A(u_j)-A(u)\right)}
  \\
  &\hphantom{\int_{\Pt^2}}\quad \times 
  \abs{f(\uDx)-f(u)} \abs{\Psi_x \omega_r \omega_{r_0} + \Psi
    \omega_r' \omega_{r_0}} \,dX\\
  &\le C\iint_{\Pt} \frac{1}{\eps} \abs{\uDx-u_j} \left(\abs{\Psi_x} +
    \Psi \frac{1}{r} \right)\,dxdt\\
  &\le \frac{CT\Dx}{\eps}\max_{t\in[0,T]} \abs{\uDx}_{B.V.(\R)}
  \left(\alpha+\frac{1}{r}\right) \\
  &\le \frac{C\Dx}{\eps r},
\end{align*}
for sufficiently small $r$ and $\alpha$. Now
\begin{align*}
  \iint_{\Pt} \abs{\eqref{eq:udx-ent-r4}}\,\bigm|_{c=u}\,dyds&\le
  M\int_{\Pt^2} \abs{\uDx-u_j} \abs{\Psi_x \omega_r \omega_{r_0} + \Psi
    \omega_r' \omega_{r_0}} \,dX \\
  &\le CT\Dx\max_{t\in[0,T]} \abs{\uDx}_{B.V.(\R)}
  \left(\alpha+\frac{1}{r}\right)\\
  &\le \frac{C\Dx}{r}.
\end{align*}
To estimate the next term we observe that $\abs{\sgn_\eps''(z)} \le C/\eps^2$,
\begin{align*}
  \iint_{\Pt} \abs{\eqref{eq:udx-ent-r5}}\,\bigm|_{c=u}\,dyds&\le
  \frac{C}{\eps^2} \iint_{\Pt} \abs{\theta_{j+1/2} - \uDx} \,dxdt\\
  &\le \frac{CT}{\eps^2} \Dx\max_{t\in[0,T]} \abs{\uDx}_{B.V.(\R)}\\
  &\le \frac{C\Dx}{\eps^2}.
\end{align*}
The next term involves $\Dp A(u_j)$ and $ A(\uDx)_x$, these can be
written
\begin{equation*}
  \Dp A(u_j)=A'(\alpha_{j+1/2})\Dp u_j\ \text{and}\ 
  A(\uDx)_x=A'(\beta_{j+1/2}) \Dp u_j,
\end{equation*}
if $x\in [x_j,x_{j+1})$. Here $a_{j+1/2}$ and $\beta_{j+1/2}$ are
between $u_j$ and $u_{j+1}$. Furthermore, since $A'\ge \eta$, we have
that
\begin{equation*}
  \abs{\Dp u_j} \le \frac{1}{\eta}\abs{\Dp A(u_j)}.
\end{equation*}
\begin{align*}
  \iint_{\Pt} \abs{\eqref{eq:udx-ent-r6}}\,\bigm|_{c=u}\,dyds &\le
  \frac{C}{\eta \eps} \iint_{\Pt} \abs{\Dp A(u_j) - A(\uDx)_x} \, dxdt \\
  &\le \frac{C}{\eta \eps} \iint_{\Pt} \abs{\alpha_{j+1/2}-\beta_{j+1/2}}
  \abs{\Dp u_j} \,dxdt\\
  &\le \frac{C}{\eps \eta^2} \iint_{\Pt} \abs{u_{j+1}-u_j} \abs{\Dp A(u_j)}
  \,dxdt\\
  &\le \frac{CT}{\eps \eta^2} \Dx \max_{t\in[0,T]} \abs{\uDx}_{B.V.(\R)}
  \\
  &\le \frac{C\Dx}{\eps\eta^2}.
\end{align*}
We continue with
\begin{align*}
    \iint_{\Pt} &\abs{\eqref{eq:udx-ent-r7}}\,\bigm|_{c=u}\,dyds\\
    &\le
    C\int_{\Pt^2} \abs{\Dp A(u_j)}\, \omega_{r_0}(x-y) \\
    &\hphantom{\le C\int_{\Pt^2}}\quad \times 
    \abs{\Psi(x_{j+1},t)\,\omega_r(x_{j+1}-y) -
      \Psi(x,t)\,\omega_r(x-y) } \,dX\\
    &\le C \iint \sum_j \int_{x_{j-1}}^{x_j}\! \int_{x}^{x_j}\!\! \abs{\Dp
      A(u_j)}\\
    &\hphantom{\le C \iint \sum_j \int_{x_{j-1}}^{x_j} \int_{x}^{x_j}
    }\!\!\!\!\!\!\!\!\!  \times 
    \abs{\Psi_x(z,t)\omega_r(z-y)+\Psi(z,t)\omega_r'(z-y)} \,dz dx dy
    dt\\
    &\le C \int \sum_j \abs{\Dp A(u_j)} \int_{x_{j-1}}^{x_j}
    \int_{x}^{x_j} \abs{\Psi_x(z,t)} + \frac{1}{r} \Psi(z,t) \,dz dx
    dt\\
    &\le CT\Dx \left(\alpha + \frac{1}{r}\right) \sup_{t\in[0,T]}
    \norm{\Dp A(u_j(t))}_{L^1(\R)} \\
    &\le \frac{C\Dx}{r},
\end{align*}
for sufficiently small $\alpha$. With similar arguments we show that
\begin{equation*}
  \iint_{\Pt} \abs{\eqref{eq:udx-ent-r8}} +
  \abs{\eqref{eq:udx-ent-r10}}\,\bigm|_{c=u} \,dyds\le 
  \frac{C\Dx}{r^2},
\end{equation*}
and 
\begin{equation*}
  \iint_{\Pt} \abs{\eqref{eq:udx-ent-r71}} +
  \abs{\eqref{eq:udx-ent-r72}}\,\bigm|_{c=u} \,dyds\le 
  \frac{C\Dx}{r\eps} + \frac{C\Dx}{r\eta}.
\end{equation*}
The term \eqref{eq:udx-ent-r9} consists of two parts. The first of
these
\begin{equation*}
  \int_{\Pt^2} \abs{\psi_\eps(\uDx,u)-\abs{\uDx-u}} \test^{\Dx}_t \,
  dX \le \frac{C\eps}{r_0\eta}\left(L+\frac{1}{\alpha}\right),
\end{equation*}
by the same arguments used to show \eqref{eq:q2est}. For the second part,
using the fact that $(\uDx)_t \in L^1(\R)$, we can show that
\begin{align*}
& \int_{\Pt^2} \abs{\uDx -u} (\test_t - (\test_j)_t)\,dX \\
& \le \int_{\Pt^2} \Bigl [ \abs{(\uDx)_t} (\psi(x_j,t) -\psi(x,t)) \omega_r(x-y) \omega_{r_0} (t-s) 
\\
& \hphantom{ \int_{\Pt^2} \Bigl[} \quad 
+ \abs{(\uDx)_t} (\omega_r(x-y) -\omega_r(x_j-y)) \psi(x,t) \omega_{r_0} (t-s) \Bigr] \,dX \\
& \le C \Dx + \frac{C \Dx}{r}.
\end{align*}
There remains the
last term, \eqref{eq:udx-ent-r11}. To estimate this, we follow the similar 
arguments to the one used in \eqref{eq:q1est}.
The result is that
\begin{equation*}
  \iint_{\Pt} \abs{\eqref{eq:udx-ent-r11}}\bigm|_{c=u}\, dyds \le
  C \frac{\eps  }{r \eta}\left(L+\frac{1}{\alpha}\right).
\end{equation*}
Now we have proved the following lemma:
\begin{lemma}
  \label{lem:estimates} Assume that $u$ and $\uDx$ are take values in
  the interval $[-K,K]$ for some positive $K$. Let $M>\max_{v\in
    [-K,K]} \abs{f'(v)}$. Then if $T\ge\tau>\nu\ge 0$ and $L-M \tau >0$, we have
  \begin{equation}
    \label{eq:estimates}
    \begin{aligned}
    \int_{-L+M\tau}^{L-M\tau} &\abs{\uDx(x,\tau)-u(x,\tau)}\,dx \\&\le 
    \int_{\R} \abs{\uDx(x,\nu)-u(x,\nu)} \,dx\\
    &\quad + C\biggl[
    r_0 + r  +\alpha
    + \Dx\left(\frac{1}{\eps^2} + \frac{1}{\eps\eta^2} +
      \frac{1}{r^2} + \frac{1}{r_0} + \frac{1}{\eps r} + \frac{1}{r} \right)\\
    &\hphantom{\quad + C\biggl[}\quad + 
      \eps \left(\frac{1}{\eta} + 
        \frac{1}{r_0\eta} + \frac{1}{r_0\eta\alpha}+ \frac{1}{r\eta\alpha} +
        \frac{1}{\eta\alpha} \right)
      \biggr].
    \end{aligned}
  \end{equation}
\end{lemma}
This follows from \eqref{eq:estimate-middle} and \eqref{eq:Rclaim},
observing that we can send $\alpha_0$ to zero.
Now we let $v(x,t)$ be the unique entropy solution of \eqref{eq:main},
where $A'(v)\ge 0$. We set $\alpha=r=r_0=\eta^{1/2}$, and assume
that $\alpha$ is sufficiently small, then
\begin{equation}
  \label{eq:vudx-estimate}
  \int_{-L+Mt}^{L-Mt} \abs{\uDx(x,t)-v(x,t)}\,dx \le
  C\left(\alpha+\frac{\eps}{\alpha^4} + \frac{\Dx}{\eps^2}\right),
\end{equation}
for some constant $C$ which is independent of the small parameters
$\alpha$, $\eps$ and $\Dx$. This follows from \eqref{eq:estimates},
\eqref{eq:diff_u_ueta} and \eqref{eq:initial_error}. Setting
$\alpha=\Dx^{1/11}$, and $\eps=\alpha^5$ proves the main theorem. 

\begin{proof}[Proof of the main corollary] To prove the corollary, we retrace
the proof of the main theorem, using that $A(u)=\eta u$. We begin by setting
\begin{equation*}
  \psi_\eps(u,c)= \eta \int_c^u \sgn_\eps(z -  c)\,dz.
\end{equation*}
In this case, the equation corresponding to \eqref{eq:udoubled} reads
\begin{multline}
  \label{eq:lin_udoubled}
  \int_{\Pt^2} \abs{u-\udx}\test_s +
  \sgn\left(u-\uDx\right)
  \left(f(u)-f\left(\uDx\right)\right) \test_y \,dX
  \\
  =\lim_{\eps \rightarrow 0}\int_{\Pt^2}\Bigl[\sgn'_\eps\left(u- \uDx \right)
  \left(u_y^2-u_y\left(\uDx\right)_x\right) \test
  \\
  -
  \abseps{u-\uDx}\left(\test_{yy}+\test_{xy}\right)
  \Bigr] \,dX.
\end{multline}
On the other hand, to obtain an expression for $\uDx$ we can proceed
as follows. In fact, in discrete set up we have the following inequality
\begin{align*}
\iint &\psi_{\eps}(u_j,c) (\test_j)_t  + Q_{\eps}(u_j,c) \Dm \test_j \,dxdt\\
& \ge \eta \iint  \Dp u_j \Dp \test_j \sgn_{\eps} (u_j -c )\,dyds 
+ \eta \iint  \partial_x [\sgn_{\eps} (\uDx -c ) ] (\uDx)_x  \test \,dxdt\\
& + \eta \iint
   \Bigl[ \Dp[\sgn_{\eps} (u_j -c ) ] \Dp u_j \test_{j+1} -
  \partial_x [\sgn_{\eps} (\uDx-c ) ] (\uDx)_x  \test \Bigr] \,dxdt\\
& = \eta \iint (\uDx)_x \test_x \sgn_{\eps}(\uDx-c) +  \sgn_{\eps}'(\uDx-c) (\uDx)_x^2 \test \, dxdt\\
& + \eta \iint (\uDx)_{xx} \sgn_{\eps} (\uDx-c) \test - \test_j \sgn_{\eps} (u_j-c) \Dm \Dp u_j  \,dxdt,
\end{align*}
where we have used $Q'_{\eps}= \psi'_{\eps} f'$ and the following equality:
\begin{align*}
0 \ge \int_{u_j}^{u_{j+1}} \psi_{\eps}^{''} (s,c) [f(u_{j+1}) -f(s)]\,ds = \int_{u_j}^{u_{j+1}}  \psi_{\eps}^{'} (s,c) & f'(s)\,ds \\
& - \psi_{\eps}'(u_j) [f(u_{j+1}) -f(u_j)]. 
\end{align*}
By taking limit as $\eps \rightarrow 0$, and choosing $c = u(y,s)$ we have
\begin{equation}
\label{eq:dis_1}
\begin{aligned}
&\int \abs{u_j-u} (\test_j)_t  +  \sgn(u_j-u) (f(u_j)-f(u)) \Dm \test_j \,dX\\
& \ge \lim_{\eps \rightarrow 0} \int \Bigl[\left((\uDx)_x^2 - u_y (\uDx)_x\right) \sgn_\eps'\left(\uDx- u \right)
     \test  - \abseps{\uDx -u}\left(\test_{xx}+\test_{xy}\right) \Bigr]\,dX.
\end{aligned}
\end{equation}
Now adding \eqref{eq:lin_udoubled} and \eqref{eq:dis_1}, using $Q(a,b)= \sgn(a-b)(f(a)-f(b))$, we get
\begin{equation}
\label{eq:dis_2}
\begin{aligned}
  \int_{\Pt^2}& \Bigl[
  \abs{\uDx-u}\left(\test_t+\test_s\right) )
  \\
  &\hphantom{ \int_{\Pt^2} \Bigl[} \quad 
  + \sgn \left( \uDx - u \right)
  \left(f\left(\uDx\right)-f(u)\right)
  \left(\test_x+\test_y\right)\notag\\
  & \hphantom{\int_{\Pt^2} \Bigl[ } 
  \qquad  +
  \abseps{\uDx -u }\left(\test_{xx}+2\test_{xy}+\test_{yy}\right)
  \Bigr]\,dX
  \\
  &\ge \int_{\Pt^2} \left( \abs{u_j-u} - \abs{\uDx-u}\right) (\test_j)_t \,dX + \int_{\Pt^2} \abs{\uDx -u} (\test_t - (\test_j)_t)\,dX\\
& + \int_{\Pt^2} Q(u_j,u) (\test_x - \Dm \test_j) \,dX + \int_{\Pt^2} \left( Q(u,\uDx) -Q(u_j,u)\right) \test_x \,dX .
\end{aligned}
\end{equation}
Observe that we can use similar arguments to the ones used in the nonlinear diffusion case for all the terms on the left side of the above inequality \eqref{eq:dis_2}. Remaining all the terms can be estimated in the following manner. We begin with
\begin{align*}
\int_{\Pt^2}  \left( \abs{u_j-u}  - \abs{\uDx-u}\right) (\test_j)_t \,dX  
 &\le \int_{\Pt^2} \abs {\uDx -u_j} \abs{ \psi_t \omega_r \omega_{r_0} + \psi \omega_r \omega'_{r_0}} \,dX \\
& \le \frac{C \alpha T \Dx}{r_0} \max_{t \in [0,T]} \abs{\uDx}_{B.V.(\R)} \\
& \le \frac{C \Dx}{r_0}.
\end{align*}
Next, we continue with
\begin{align*}
\int_{\Pt^2}& Q(u_j,u)  (\test_x - \Dm \test_j) \,dX  = - \int_{\Pt^2} Q(u_j,u)_x \test - \Dp Q(u_j,u) \test_j \,dX\\
& = \int_{\Pt} \int_0^T \sum_j \int_{x_j}^{x_{j+1}} \left( \Dx \test \delta_{\lbrace x_{j+1/2} \rbrace } \Dp Q(u_j,u) - \Dp Q(u_j,u) \test_j \right) \,dX \\
& \le \frac{C \Dx}{r}
\end{align*}
In a similar way, we also can show that
\begin{align*}
\int_{\Pt^2} \left( Q(u,\uDx) -Q(u_j,u)\right) \test_x \,dX \le \frac{C \Dx}{r}
\end{align*}
Finally, we end up with
\begin{equation}
    \label{eq:linearestimates}
    \begin{aligned}
    \int_{-L+M\tau}^{L-M\tau} &\abs{\uDx(x,\tau)-u(x,\tau)}\,dx \\&\le 
    \int_{\R} \abs{\uDx(x,\nu)-u(x,\nu)} \,dx 
     + C \left(
    r_0 + r  +\alpha
    + \frac{\Dx}{r} + \frac{\Dx}{r_0} \right).
    \end{aligned}
  \end{equation}
Now setting $\alpha=r=r_0$ yields the estimate
\begin{equation*}
   \int_{-L+Mt}^{L-Mt} \abs{\uDx(x,t)-v(x,t)}\,dx \le C\left(\alpha + 
\frac{\Dx}{\alpha}\right).
\end{equation*}
Setting $\alpha=\Dx^{1/2}$ proves the corollary.
\end{proof}

\section{A numerical test}\label{sec:numtest}
In order to test the unlikely optimality of the convergence rate of
our main theorem, we compute the numerical convergence rate of an
example. Consider the following initial value problem
\begin{equation}
  \label{eq:num_initial}
  \begin{cases}
    u_t = A(u)_{xx}\
    &\text{for $t>0$ and $x\in (-\pi/2,\pi)$,}\quad 
    A(u)= \frac{1}{2} \left(\max\seq{u,0}\right)^2\\
    u(x,0)=\sin(x), & x\in [-\pi/2,\pi],
  \end{cases}
\end{equation}
supplemented with the boundary conditions 
\begin{equation*}
  \partial_x A(u(t,x))=0\quad \text{for $t>0$ and $x=-\pi/2$, $x=\pi$.}
\end{equation*}
In order to simplify matters, we have chosen an example without the
convective term, nevertheless a discontinuity will form in the
solution due to the degeneracy of $A$. This discontinuity will take
the form of a boundary on which $u=0$,  moving to the left. 

We have used the Euler method to integrate the system of ordinary
differential equations \eqref{eq:scheme}, resulting in the update
formula
\begin{equation*}
  u_j((n+1)\Dt)=u_j(n\Dt) + \Dt \Dm\Dp A_j(n\Dt).
\end{equation*}
For this to be linearly stable, the time-step $\Dt$ must obey
the restrictive CFL condition $\Dt \le 0.5 \max A'(u) \Dx^2$. 

\begin{figure}[tph]
  \centering
  \begin{tabular}[h!]{rl}
    \includegraphics[width=0.5\linewidth]{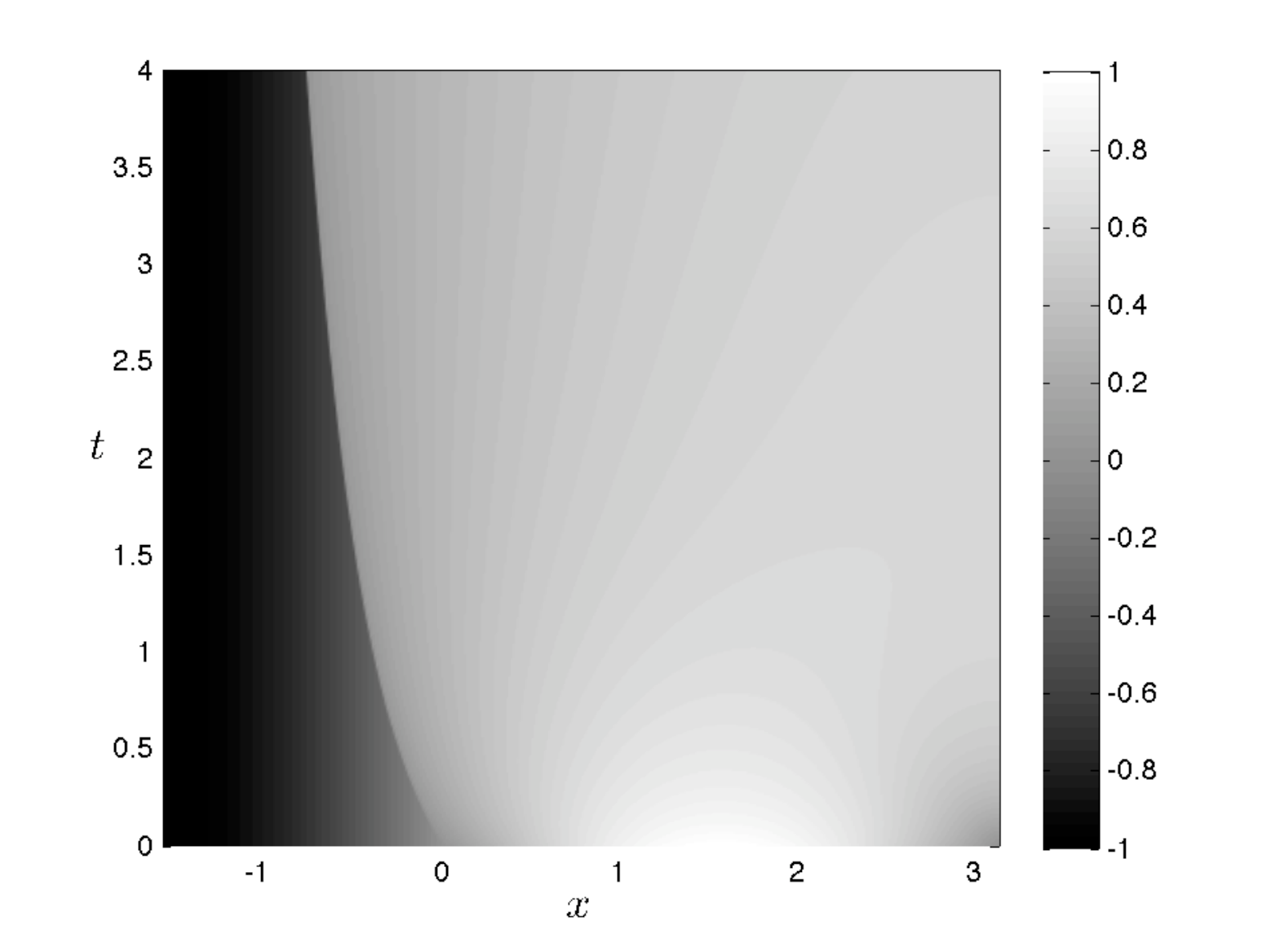}&
    \includegraphics[width=0.5\linewidth]{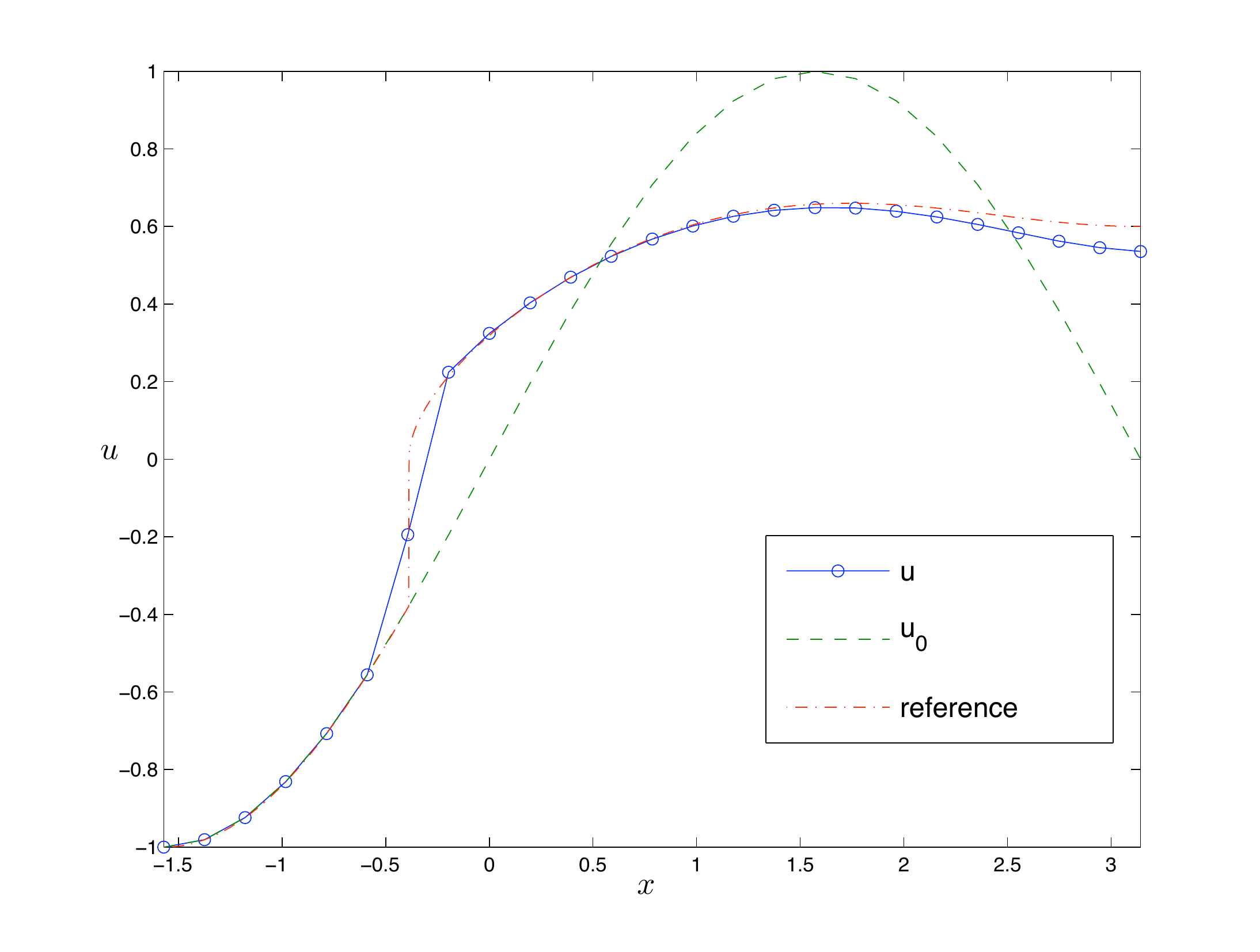}
  \end{tabular}
  \caption{An approximate solution to \eqref{eq:num_initial} using
    $400$ grid points. Left: $u$ in the $(x,t)$ plane for $t\in
    [0,4]$. Right: an approximation to $u(1,x)$ using $25$ grid
    points, a reference solution computed using $4000$ grid points and
    the initial data. }
  \label{fig:n1}
\end{figure}

In
Figure~\ref{fig:n1} we show the solution in the $(x,t)$ plane and a
snapshot of $u$ at $t=1$, for an approximation using $400$ grid points
in the interval $(-\pi/2,\pi)$. 

Finally we computed approximate errors, defined by 
\begin{equation*}
  100\frac{\norm{u_{\Dx} - u_{\mathrm{ref}}}_{L^1}}{\norm{u_{\Dx}}_{L^1}}
\end{equation*}
for a reference solution computed by using our scheme with $4000$ grid
points. These errors were computed at $t=1$. The result of this is
shown in Table~\ref{tab:1}.
\begin{table}[thp]
  \centering
  \begin{tabular}[h!]{c|rrrrrr}
    $N$&\multicolumn{1}{c}{$25$}
    &\multicolumn{1}{c}{$50$}
    &\multicolumn{1}{c}{$100$}
    &\multicolumn{1}{c}{$200$}
    &\multicolumn{1}{c}{$400$}
    &\multicolumn{1}{c}{$800$}
    \\
    \hline
    error&3.62&1.55&0.82&0.40&0.18&0.07\\
    rate &\multicolumn{1}{c}{-}&
    1.22&0.92&1.02&1.11&1.42
  \end{tabular}
  \caption{Numerical errors and convergence rate for the intial value
    problem \eqref{eq:num_initial}. $N$ is the number of grid points.}
  \label{tab:1}
\end{table}
Looking at this table, it seems the numerical rate of convergence
is $1$. Thus very far from the rather pessimistic lower bound proved
in this paper!

\bibliographystyle{abbrv}
\bibliography{Ujjwal}  

\end{document}